\documentclass[11pt]{article}
\usepackage{amssymb, amsmath, theorem, epsfig, latexsym}

\numberwithin{equation}{section}

\theoremstyle{plain}
\theorembodyfont{\itshape}
\newtheorem{theorem}{Theorem}[section]
\newtheorem{proposition}[theorem]{Proposition}
\newtheorem{lemma}[theorem]{Lemma}

\theorembodyfont{\rmfamily}
\newtheorem{definition}[theorem]{Definition}

\newtheorem{example}[theorem]{Example}
\newtheorem{remark}[theorem]{Remark}
\newenvironment{proof}{{\noindent \textbf{Proof}\,\,}}{\hspace*{\fill}$\Box$\medskip}

\title{\textbf{Stokes matrices of hypergeometric integrals}}
\author{Alexey Glutsyuk\thanks{Laboratoire
J.-V.Poncelet (UMI 2615 du CNRS et l'Universit\'e Ind\'ependante
de Moscou). Permanent address: CNRS, Unit\'e de Math\'ematiques
Pures et Appliqu\'ees, M.R., \'Ecole Normale Sup\'erieure de Lyon,
46 all\'ee d'Italie, 69364 Lyon 07, France.  \newline Email:
aglutsyu@umpa.ens-lyon.fr }, \ Christophe
Sabot\thanks{Universit\'e de Lyon, Universit\'e Lyon 1, Institut
Camille Jordan, CNRS UMR 5208, 43 bd du 11 nov. 1918, 69622
Villeurbanne Cedex.
\newline Email: sabot@math.univ-lyon1.fr }}
\begin{document}
\maketitle

\def\Rm{\refstepcounter{remarques}{\bf Remark \theremarques} }
\def\wt{\widetilde}
\def\ali{\hfill\break}
\def\w{\omega}
\def\W{\Omega}
\def\unN{\{1\ldots ,N\}}
\def\rrr{{\cal R}}
\def\demi{{1\over 2}}
\def\ddd{{\cal D}}
\def\vvv{{\cal V}}
\def\nn{{<n>}}
\def\ND{{ND}}
\def\infi{{<\infty>}}
\def\intX{{\stackrel{\circ}{X}}}
\def\npn{{<n+p>,<n>}}
\def\unN{{\{1,\ldots ,N\}}}
\def\ccc{{\cal C}}
\def\j1p0{{j_1^0,\ldots ,j_p^0}}
\def\intccc{{\stackrel{\circ}{{\cal C}}}}
\def\jp1{{j_p,\ldots ,j_1}}
\def\j1p{{j_1,\ldots ,j_p}}
\def\dn{{{\partial_n}}}
\def\Ppk{{{I\!\!\! P}^k}}
\def\oomega{{\overline\omega}}
\def\Z{{{\Bbb Z}}}
\def\P{{{\Bbb P}}}
\def\N{{{\Bbb N}}}
\def\R{{{\Bbb R}}}
\def\C{{{\Bbb C}}}
\def\E{{{\Bbb E}}}
\def\Fn{{F^{(n)}}}
\def\Fnmoinsun{{F^{(n-1)}}}
\def\Finfi{{F^{(\infty)}}}
\def\F1{{{F^{(1)}}}}
\def\En{{E^{(n)}}}
\def\Gn{{G^{(n)}}}
\def\Yn{Y^{(n)}}
\def\Ynk{{Y^{(n)}_k}}
\def\eee{{{\cal E}}}
\def\eeeun{{{\cal E}^{(1)}}}
\def\supp{{{\hbox{supp}}}}
\def\ND{{{ND}}}
\def\symG{{{\hbox{Sym}^G}}}
\def\Fun{{{F_{<1>}}}}
\def\oeta{{{\overline{\eta}}}}
\def\aaa{{{\cal A}}}
\def\un{{<1>}}
\def\mult{{\hbox{mult}}}
\def\ppp{{{\mathcal P}}}
\def\aaa{{{\mathcal A}}}
\def\lll{{{\mathcal L}}}
\def\clos{{{\hbox{closure}}}}
\def\iii{{{\mathcal I}}}
\def\cad{c'est-\`a-dire }
\def\niii{{{\not\iii}}}
\def\Sl{{{\hbox{Sl}}}}
\def\iff{if and only if}
\def\card{{\hbox{card}}}
\def\fff{{\mathcal F}}
\def\WFP{{(\Omega, \fff, \P)}}
\def\var{\hbox{var}}
\def\indic{{\bf 1}}
\def\equalinlaw{{\stackrel{\hbox{law}}=}}
\def\hhh{{\mathcal H}}
\def\ccc{{\mathcal C}}
\def\bbb{{\mathcal B}}
\def\diag{{\hbox{diag}}}
\def\xxx{{\mathcal X}}
\def\sgn{{\hbox{sgn}}}
\def\ddd{{\mathcal D}}
\def\re{\operatorname{Re}}
\def\im{\operatorname{Im}}

{\footnotesize \noindent{\slshape\bfseries Abstract.}} In this
work we compute the Stokes matrices of the ordinary differential
equation satisfied by the hypergeometric integrals associated to 
an arrangement of hyperplanes in generic position. This generalizes
the computation done by Ramis and Duval for confluent
hypergeometric functions, which correspond to the arrangement of
two points on the line. The proof is based on an explicit
description of a base of canonical solutions as integrals on the
cones of the arrangement, and combinatorial relations between
integrals on cones and on domains.

\section{Introduction and main result}

The computation of the Stokes matrix of an ordinary differential
equation with an irregular singular point is in general a
difficult problem. In \cite{Ramis1} and \cite{Ramis2}, Ramis and
Duval considered the case of confluent hypergeometric functions,
and computed the associated Stokes matrices. In this paper, we
consider a natural generalization: we consider an arrangement of
hyperplans in generic position and the hypergeometric integrals
with an exponential term of the form $e^{-\lambda f_0}$ where
$f_0$ is an extra linear form. Differentiating in $\lambda$ leads
to a differential equation satisfied by these integrals, with a
regular singular point at 0 and an irregular singular point at
infinity. The case of \cite{Ramis1,Ramis2} is the case of the
arrangement of two points on the line. The purpose of this paper
is to compute explicitly the stokes matrices of this equation. A
differential equation of this type appears in the analysis of a
probabilistic model of random environments (\cite{sabot}), which
was one of the motivation of this work.

Let $f_1, \cdots ,f_N$ be $N$ affine forms on $\R^k$, $N\ge k$,
and set
$$
H_j=\ker f_j.
$$
We assume that the hyperplanes $H_1, \ldots ,H_N$ are in generic
position (all of them are distinct, any $k$ planes intersect at a single
point and the intersection of any $k+1$ planes is empty). We denote by
$$
l_j(z)=f_j(z)-f_j(0),
$$
the linear form directing $f_j$. We associate a positive weight
$\alpha_j$ to each hyperplane $H_j$, and for any subset $U\subset
\unN$ we set
\begin{eqnarray}\label{aU}
 \alpha_U=\sum_{j\in U} \alpha_j .
\end{eqnarray}
The couple $(\R^k,(H_j)_{j=1,\cdots, N})$ defines an arrangement
of hyperplanes. To any collection of $k$ hyperplanes
$H_{j_1}, \cdots ,H_{j_k}$, $j_l\neq j_r$ for $l\neq r$, we associate the unique
vertex of the arrangement
\begin{eqnarray}\label{vertex}
X=H_{j_1}\cap \cdots \cap H_{j_k}.
\end{eqnarray}
Depending on the context we will consider a vertex as a subset of
$\unN$ with $k$ elements (i.e. in (\ref{vertex}), $X=\{j_1, \ldots
,j_k\}$) or a point of $\R^k$ (as in formula (\ref{vertex})). We
denote by $\xxx$ the set of vertices of the arrangement. To any
vertex $X=\{j_1, \ldots j_k\}$ we associate the differential form
of maximal degree
$$
\w_X={\frac{df_{j_1}}{f_{j_1}}}\wedge \cdots \wedge {\frac{df_{j_k}}
{f_{j_k}}},
$$
where the elements of $X$ are ordered so that the form
$df_{j_1}\wedge \cdots \wedge df_{j_k}$ is positively oriented
(for an arbitrary fixed orientation of the vector space $\R^k$).

\renewcommand\Re{\operatorname{Re}}

A connected component $\Delta$ of $\R^k\setminus \cup_{j=1}^N
H_j$ is called an {\it arrangement domain}. We denote by $\ddd$ the set of the
arrangement domains. Let $f_0$ be a linear form on $\R^k$ in general
position with respect to $(f_1, \ldots ,f_N)$ (i.e., $f_0$ takes distinct values on the
vertices of the arrangement and is nonconstant on each intersection line
of $k-1$- ple of hyperplanes). We denote by
$\ddd^+$  the set of the arrangement domains on which the form $f_0$ is bounded
from below. Since the arrangement is generic, it follows that the
domains of $\ddd^+$ are the bounded domains or the unbounded
domains $\Delta$ such that there exist some constants $A\in\mathbb R$ and
$B>0$ such that $f_0(x)\ge A+B \|x\|$ on $\Delta$. To any domain
$\Delta$ of $\ddd^+$ and any vertex $X$, we associate the
integral
\begin{equation}
I_{\Delta,X} (\lambda)=\int_\Delta e^{-\lambda f_0}\Omega_X, \
\Omega_X=\left(\prod_{j=1}^N \vert f_j\vert^{\alpha_j}\right) \w_X,
\label{int}\end{equation}
for $\Re(\lambda)>0$.

Now we need to describe the edges of dimension 1 of the
arrangement: to any subset $U=\{j_1, \ldots ,j_{k-1}\}\subset X$ we
associate the edge of the arrangement
$$
L_U=\cap_{j\in U} H_j,
$$
which is a line in $\R^k$. Let $e_U$ be the unique vector
directing $L_U$, i.e. such that $L_U=X+\R e_U$, and normalized so
that
\begin{eqnarray}\label{eU}
f_0(e_U)=1.
\end{eqnarray}

 The general theory of hypergeometric integrals tells that
these integrals are solutions of a differential equation. In our
case, we can show (for the convenience of the reader, we give a proof
of this result at the end of the paper) that for any domain
$\Delta$ in $\ddd^+$, the vector
$$
I_\Delta(\lambda)=(I_{\Delta,X}(\lambda))_{X\in \xxx}
$$
satisfies the following ordinary differential equation
\begin{eqnarray}\label{ODE}
I'=-(\aaa+\frac1{\lambda} \bbb) I,
\end{eqnarray}
where $\aaa$ is the diagonal matrix with diagonal terms
$$
\aaa_{X,X}=f_0(X).
$$
The matrix $\bbb$ is given by
$$
\bbb_{X,X}=\alpha_X
$$
on the diagonal and
$$
\bbb_{X,Y}=0,
$$
if the vertices $X$, $Y$ are distinct and do not lie in one and the same edge (or
equivalently, $\vert X\cap Y\vert<k-1$).

Finally, if $\vert X\cap Y\vert =k-1$, we set $\{j\}=X\setminus
Y$, $\{r\}=Y\setminus X$, $U=X\cap Y$,
$$
\bbb_{X,Y}= \epsilon(j,r,U)\alpha_r,
$$
where $\epsilon(j,r,U)$ depends on the relative orientation of
$f_j$ and $f_r$ on the edge $L_U$:
$$
\epsilon(j,r,U)= \sgn(l_j(e_U)l_r(e_U)).
$$

There is a natural bijection between the set of vertices $\xxx$
and the domain set  $\ddd^+$: to each domain $\Delta\in
\ddd^+$ we associate the unique vertex $X(\Delta)\in
\partial\Delta$ that minimizes $f_0$ on $\Delta$; the inverse of
this application associates to any vertex $X$ the unique domain
$\Delta_X$ containing $X$ in its boundary and on which
$f_0-f_0(X)>0$. The Wronskian of the solutions
$(I_\Delta(\lambda))_{\Delta\in \ddd^+}$ has been explicitely
computed in the works by A.N.Varchenko
\cite{Varchenko1, Varchenko10}, in his joint work with Y.Markov
and V.Tarasov \cite{Varchenko2}, and in the joint work by A.Douai and
 H.Terao \cite{terao}. This Wronskian is nonzero.
Hence, the functions $(I_\Delta(\lambda))_{\Delta\in
\ddd^+}$ form a basis of solutions of the differential system
(\ref{ODE}) on the set $\{\Re(\lambda)>0\}$. The differential
equation (\ref{ODE}) admits a regular singular point at
$\lambda=0$ and an irregular singular point at $\lambda=\infty$.
The question we address in this paper is the explicit computation
of the Stokes matrices of this differential equation. J.-P.Ramis \cite{Ramis1} and A.Duval \cite{Ramis2} computed the Stokes matrices
of some confluent hypergeometric integrals, which corresponds to a
particular case of our differential equations (cf. Example \ref{ex1}).

\def\var{\varepsilon}

 The general theory (see \cite{il, kh}) says that there is a unique formal linear invertible
 change of space variables at infinity
 that transforms (\ref{ODE}) to its formal normal form:
 \begin{equation}Y'=-(\aaa+\frac1{\lambda}diag(\bbb))Y,\label{nform}\end{equation}
 where $\diag(\bbb)$ is the diagonal matrix formed by the diagonal
terms of $\bbb$ (i.e. $\bbb_{X,X}=\alpha_X$). The previous formal change
is given by a formal Laurent nonpositive power
 series in $\lambda$ (with matrix coefficients; the free term is unit) that does not
converge in general. On the other hand, on each sector $S_{\pm}\subset\mathbb C$
defined below there exists a unique holomorphic variable change
(called {\it sectorial normalization}) transforming (\ref{ODE})
to (\ref{nform}) for which the previous normalizing series is its asymptotic Laurent
series at infinity. The latter statement holds true for the following sectors, see Fig.1a:
\begin{equation}S_{\pm}=\{\var-\frac{\pi}2<\pm\arg\lambda<\frac{3\pi}2-\var\}; \
\text{with arbitrarily fixed} \ \var, \ 0<\var<\frac{\pi}2.\label{sectors}
\end{equation}
\begin{definition} \label{defbas} The {\it canonical solution base} of (\ref{nform}) is the
base of its solutions given by a diagonal fundamental matrix.
The {\it canonical sectorial solution base} of (\ref{ODE}) in $S_{\pm}$
is its pullback under the corresponding sectorial normalization.
\end{definition}

The canonical solution bases are uniquely defined up to multiplication of
the base solutions by constants. We normalize them as follows. Let
$$V\to\mathbb C^*=\mathbb C\setminus0$$
 be the universal
cover over $\mathbb C^*$. We lift both equations (\ref{ODE}) and
(\ref{nform}) and the sectorial normalizations  to $V$.
Take a holomorphic branch on $V$ of the diagonal
fundamental solution matrix $W$ of the formal normal form (\ref{nform}). Fix
connected components $S_0,S_1,S_2\subset V$ of the covering projection
preimages  of
$S_+$, $S_-$ and $S_+$ respectively that are ordered clockwise so that
\begin{equation}S_{01}=S_0\cap S_1\neq\emptyset, \
S_{12}=S_1\cap S_2\neq\emptyset , \ \text{see Fig.1b}.\label{nsectors}
\end{equation}
\begin{definition} \label{defnbas} The {\it normalized tuple of canonical sectorial
solution bases} of
equation (\ref{ODE})  in $S_j$, $j=0,1,2$, consists of the pullbacks of the
previous holomorphic fundamental matrix $W$ under the corresponding
sectorial normalizations of (\ref{ODE}). Then for any $j=0,1$ the pair of the
previous solution bases in $S_j$ and $S_{j+1}$ is called a {\it normalized base pair}.
\end{definition}
\begin{figure}[ht]
  \begin{center}
   \epsfig{file=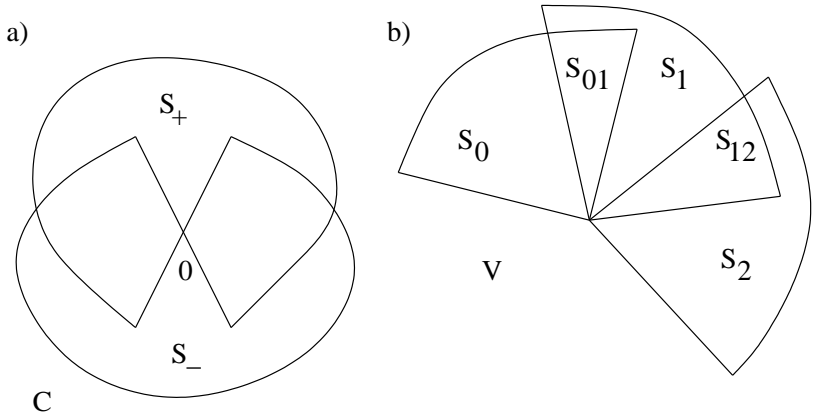}
    \caption{The sectors $S_{\pm}$, $S_0$, $S_1$, $S_2$}
    \label{fig:1}
  \end{center}
\end{figure}
\begin{remark} \label{rnbas} A normalized base tuple (pair) is uniquely defined up to
multiplication of the base functions by constants (independent on the sector).
\end{remark}

Denote the previous normalized sectorial solution bases in $S_j$
(more precisely, their fundamental matrices) by $Z_j(\lambda)$, $j=0,1,2$.
The transitions between them  in the intersections $S_{01}$, $S_{12}$ of their
definition domains
are given by constant matrices $C_0$, $C_1$ called {\it Stokes matrices}:
\begin{equation}Z_1(\lambda)=Z_0(\lambda)C_0 \ \text{in} \ S_{01}, \
Z_2(\lambda)=Z_1(\lambda)C_1 \ \text{in} \ S_{12}.\label{stokes}
\end{equation}
\begin{remark} The Stokes matrices are uniquely defined up to  simultaneous
conjugation by one and the same diagonal matrix.
\end{remark}

In the present paper we find explicitly the above canonical sectorial solution bases
(Proposition \ref{canbas} in the next Section) and calculate
the corresponding Stokes matrices (the next Theorem).

We order all the vertices $X$ of the hyperplane arrangement by the corresponding
values $f_0(X)$ of the linear function $f_0$ (which are distinct by definition).
The sectorial solution bases given by Proposition \ref{canbas} in the sectors
$S_{\pm}$ are numerated by the vertices $X$. Their $X'$- components are given by
the integrals $I_{X,X'}^{\pm}$ over appropriate cones based at $X$
of the (appropriately extended) forms $e^{-\lambda f_0}\Omega_{X'}$.

To describe the Stokes matrices, we need to introduce some
notations. Let $X$ be a vertex, we denote by $\ccc_X^+$ the unique (open)
cone defined by the hyperplanes $(H_j)_{j\in X}$ on which
$f_0-f_0(X)$ is positive. Similarly, the cone $\ccc_X^-$ is the
unique cone defined by the hyperplanes $(H_j)_{j\in X}$ on which
$f_0-f_0(X)$ is negative.

\begin{definition} \label{dexcep} A pair $(X,X')$ of distinct vertices $X,X'\in\xxx$
is said to be {\it positive exceptional}, if either $X'\notin\overline\ccc^+_X$, or
$X'\in\ccc^+_X$ and there exists an arrangement hyperplane through $X'$ that does not
separate the domains $\Delta_{X'}$ and $\Delta_X$ (see Fig. 2). The latter hyperplane
is then also called {\it exceptional}. A pair $(X,X')$ is said to be {\it negative exceptional},
if it is positive exceptional with respect to the arrangement equipped with the new
linear function $\wt f_0=-f_0$.
\end{definition}
\begin{figure}[ht]
  \begin{center}
   \epsfig{file=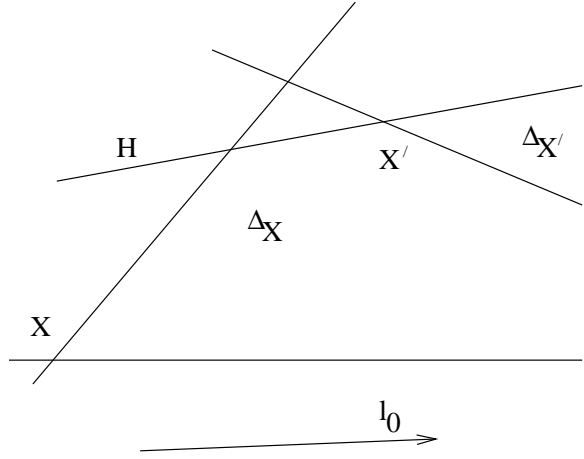}
    \caption{A positive exceptional pair $(X,X')$: the  line $H$ is exceptional}
    \label{fig:2}
  \end{center}
\end{figure}
\begin{theorem} \label{tstokes}
Consider the normalized tuple of canonical sectorial solution bases in $S_0$, $S_1$,
$S_2$ (numerated by the vertices $X\in\xxx$)
given by Proposition \ref{canbas}. The corresponding Stokes matrices
$C_j=(C_j(X',X))_{X',X\in\xxx}$, $j=0,1$,
see (\ref{stokes}), are given by the following formulas:
$$C_0(X,X)=C_1(X,X)=1,$$
\begin{equation}
C_0(X',X)=\begin{cases} & 0, \ \text{if} \ \text{the pair} \ (X,X') \
\text{is positive exceptional; \ otherwise}\\
& (-1)^{|B|+|X'\setminus X|}e^{\pi i(\alpha_{B}-\alpha_A)}
\prod_{j\in X'\setminus X}(2i\sin\pi\alpha_j);
\end{cases}\label{co}\end{equation}
\begin{equation}
C_1(X',X)=
\begin{cases} &
0 \ \text{if} \ \text{the pair} \ (X,X') \ \text{is negative exceptional; otherwise}\\
&  (-1)^{|B|+|X'\setminus X|}e^{\pi i(\alpha_{X}-\alpha_{X'}+\alpha_{B}-\alpha_A)}
\prod_{j\in X'\setminus X}(2i\sin\pi\alpha_j),
\end{cases}
\label{c1}\end{equation}
where
\begin{equation}A=\{ j \ | \ H_j \ \text{separates (strictly)} \ X \ \text{from} \ X'\},
 \label{ab}\end{equation}
 $$
B=\{ j \ | \ H_j \ \text{contains} \ X, X' \ \text{and separates the cone} \
\ccc_X^{+} \ \text{from} \ \ccc_{X'}^{+}\}.
$$
\end{theorem}

\begin{remark} The above set $B$ coincides with the set defined in a similar way
but with the upper index "$+$" of the cones replaced by "$-$". Indeed, any given
hyperplane $H$ through $X$ and $X'$ that separates the cones
$\ccc_X^+$ and $\ccc_{X'}^+$ also separates $\ccc_X^-$ from $\ccc_{X'}^-$
and vice versa. This follows from the fact that the central symmetry with respect
to $X$ ($X'$) sends $\ccc^+_X$ to $\ccc^-_X$ (respectively, $\ccc^+_{X'}$ to
$\ccc^-_{X'}$) and changes the side of the cone under consideration with
respect to $H$.
\end{remark}

\begin{example} \label{ex1} Let $k=1$, and $X_1<\cdots <X_N$ be $N$ points on
the real line, and set
$$
f_i(z)=z-X_j, \;\; \; j=1, \ldots ,N, \;\; z\in \R,
$$
$$
f_0(z)=z.
$$
The matrix $\aaa$ is the diagonal matrix
$$
\aaa=\left( \begin{array}{ccc} X_1&& \Large{0}\\ & \ddots &
\\ \Large{0}& & X_N \end{array}\right),
$$
and
$$
\bbb=\left( \begin{array}{cccc} \alpha_1 & \alpha_2 & \cdots & \alpha_N \\
\alpha_1 & \alpha_2 & \cdots & \alpha_N
\\
\vdots & \vdots & \vdots &\vdots
\\
\alpha_1 & \alpha_2 & \cdots & \alpha_N
\end{array}\right).
$$
The Stokes matrices are
$$
C_0= \left( \begin{array}{ccccc}1& 0 &\cdots & 0
\\
{-2i\sin\pi\alpha_2}
& 1& \cdots & 0
\\
 \vdots &\vdots &\cdots &\vdots
\\
{-2ie^{-\pi i\sum_{s=2}^{N-1}\alpha_s}\sin\pi\alpha_N} &
{-2ie^{-\pi i \sum_{s=3}^{N-1}\alpha_s}\sin\pi\alpha_N}  &\cdots &1
\end{array}\right),
$$
$$C_1= \left( \begin{array}{ccccc}1& {-2ie^{\pi i(\alpha_2-\alpha_1)}\sin\pi\alpha_1}
 &\cdots &
{-2ie^{\pi i(\alpha_N-\sum_{j=1}^{N-1}\alpha_j)}\sin\pi\alpha_1}
\\
0 & 1&  \cdots &
{-2ie^{\pi i(\alpha_N-\sum_{j=2}^{N-1}\alpha_j)}\sin\pi\alpha_2} \\

\\ \vdots &\vdots &\cdots &\vdots\\
0 & 0& \cdots&1
\end{array}\right).
$$
The case where $N=2$ and $X_1=0$, $X_2=1$ corresponds to the usual
confluent hypergeometric case, which has been considered in
\cite{Ramis1, Ramis2}.
\end{example}

\begin{example} \label{ex2} Let $k=2$ and for $z=(x,y)$
$$
f_1 (z)=x, \;\; f_2(z)=y, \;\; f_3(z)=x+y-1,
$$
$$
f_0(z)=ax+by,
$$
with $a>0$, $b>0$, $a>b$. The vertices of the arrangement are
$$
X_{1}:= (0,0), \;\; X_{2}:=(0,1), \;\; X_{3}:=(1,0).
$$
We have
$$
f_0(X_{1})=0<f_0(X_{2})=b< f_0(X_{3})=a,
$$
$$X_1=\{1,2\}, \ X_2=\{1,3\}, \ X_3=\{2,3\},$$
$$A=\emptyset \ \text{for each pair} \ (X',X),$$
$$B=\{ H_2\} \ \text{if} \ \{ X',X\}=\{ X_1,X_3\}, \ B=\emptyset \
\text{otherwise},$$
$$ C_0= \left( \begin{array}{ccc}1& 0&0
\\
{-2i\sin\pi\alpha_3} & 1 & 0
\\
 {2ie^{\pi i\alpha_2}\sin\pi\alpha_3}&{-2i\sin\pi\alpha_2} &1
\end{array}\right),
$$
$$
C_1= \left( \begin{array}{ccc}1& {-2ie^{i\pi(\alpha_3-\alpha_2)}\sin\pi\alpha_2} & {2ie^{\pi i(\alpha_2+\alpha_3-\alpha_1)}\sin\pi\alpha_1}
\\
0& 1 & {-2ie^{i\pi(\alpha_2-\alpha_1)}\sin\pi\alpha_1}
\\
0 & 0 &1
\end{array}\right).
$$
\end{example}


\section{Canonical solutions at infinity. The plan of the proof of Theorem \ref{tstokes}}
\subsection{Canonical solutions}
Let $X$ be a vertex and $\rho\in\C$, $|\rho|=1$. We denote by
$\ccc^\rho_X\subset \C^k$ the cone based at $X$ and defined by
\begin{equation}
\ccc_X^\rho=\{ z=X+\rho(\sum_{j\in X} a_j
e_{X\setminus\{j\}}),\;\; a_j\in \R_+\},\label{crho}\end{equation}
where $e_{X\setminus\{j\}}$ has been defined in (\ref{eU}).

 \begin{remark} When $\rho=\pm1$, one has $\ccc_X^{\rho}=\ccc_X^{\pm}$
(the cones $\ccc_X^{\pm}$ were defined in the Introduction, just before
Theorem \ref{tstokes}). For any vertex $X$ one has
$\Delta_X\subset\ccc^+_X=\ccc^1_X$.
\end{remark}

 For any affine subspace $H\subset\mathbb R^k$ denote
 $\mathbb CH\subset\mathbb C^k$ its complexification.

 \begin{remark} \label{conconf} For any
$j\in X$, the intersection $\mathbb CH_j\cap \overline{\ccc_X^\rho}$ is a face
of the cone $\ccc^\rho_X$. For any $\rho\not\in \R$ and $l\notin X$ one has
$$\overline{\ccc_X^\rho}\cap \mathbb C H_l=\emptyset.$$
Without loss of generality we prove this statement assuming that $X=0$
(translating the coordinates). Suppose the  contrary: there exist a $\rho\in\mathbb R$
and a $l\notin X$ such that there exists a point $x_0\in\overline{\ccc_0^{\rho}}\cap
\mathbb CH_l$. By definition, $0=X\notin\mathbb CH_l$, since
$0\notin H_l$ ($l\notin X$ by assumption). In particular, $x_0\neq 0$. One has
\begin{equation} x_0=\rho v,  \ v\in\overline{\ccc^+_0}\setminus0, \ x_1=
\re x_0=(\re\rho)v\in\re(\mathbb CH_l)=H_l,\label{x0}\end{equation}
\begin{equation}x_2=\im x_0=(\im\rho)v\in\im(\mathbb CH_l)=H_l', \label{x2}
\end{equation}
where $H_l'$ is the real hyperplane  through 0 parallel to $H_l$. One has
$x_2\neq0$, since $x_0\neq0$ and $\im\rho\neq0$ ($\rho\notin\mathbb R$ by
assumption). The vector $x_1$ lies in $H_l'$, since it is proportional to
$x_2\in H_l'\setminus0$. Therefore, $x_1$ lies simultaneously in two disjoint hyperplanes
$H_l$ and $H_l'$, - a contradiction.
\end{remark}

For any $\rho\not\in \R$, we consider the integral
\begin{eqnarray}\label{IX}
I^{\rho}_{X,X'} (\lambda)=\int_{\ccc_X^\rho} e^{-\lambda f_0}\Omega_{X'},
\ \Omega_{X'}=\left( \prod_{j=1}^N |f_j|^{\alpha_j}\right) \w_{X'},
\end{eqnarray}
\def\cxr{\ccc_x^\rho}
\def\Im{\operatorname{Im}}
where the determination of the 1- form $\Omega_{X'}$ is chosen as follows.
Take a simply connected domain $\widetilde{\mathcal D}\subset\mathbb C^k\setminus
\cup_j\mathbb CH_j$ containing the union of the cones $\ccc^{\rho}_X$, $\im\rho<0$.
(The latter cones are simply connected, as is their union, and disjoint from
the complex hyperplanes $\mathbb CH_j$ (see the previous Remark). Hence, the
previous domain
$\widetilde{\mathcal D}$ exists.) Take the standard real branch of $\Omega_{X'}$ on the real
domain $\Delta_X\subset\ccc^1_X=\ccc^+_X$. The domain $\Delta_X$ lies in
$\mathbb R^n\setminus\cup_jH_j$ and is adjacent to the previous union of cones. Take
the immediate analytic extension of the real branch  $\Omega_{X'}|_{\Delta_X}$ to
$\widetilde{\mathcal D}$.

\begin{remark}
The integral (\ref{IX}) is well-defined whenever $\lambda$ is such
that $\Re ( \lambda \rho)>0$. Moreover, for any $\lambda\not\in
i\R_-$, the integral does not depend on $\rho$ such that $\Im
\rho<0$ and $\Re(\rho\lambda)>0$ (when $\lambda \in i\R_-$,
there is no such $\rho$).
\end{remark}

We denote by
$$
I^{+}_{X,X'}(\lambda)
$$
the common value of $I^{\rho}_{X,X'} (\lambda)$ for $\Im(\rho)<0$
and $\Re(\rho\lambda)>0$. The function $I^{+}_{X,X'}(\lambda)$ is
analytic on $\C\setminus i\R_-$. Similarly, we denote by
$I^{-}_{X,X'}(\lambda)$ the common value of
$I^{\rho}_{X,X'}(\lambda)$ for $\Im (\rho) >0$ and
$\Re(\rho\lambda)>0$. The function $I^{-}_{X,X'}(\lambda)$ is
well-defined and analytic on $\C\setminus i\R_+$.

We denote by
$I^{\pm}_{X}(\lambda)$ the vector
$$
I^{\pm}_{X}(\lambda)=(I^{\pm}_{X,X'}(\lambda))_{X'\in \xxx}, \
\Omega=(\Omega_{X'})_{X'\in\xxx}.
$$

\def\var{\varepsilon}

\begin{proposition} \label{canbas}
Let $S_{\pm}$, $V$, $S_0,S_1,S_2\subset V$ be as in (\ref{nsectors}).
The vector functions $I^{\pm}_X(\lambda)$
(corresponding to all the vertices $X$) are solutions of
(\ref{ODE}) and form a canonical sectorial solution basis in the corresponding
sector $S_{\pm}$ (see (\ref{sectors})). The liftings to $S_0$, $S_1$, $S_2$ of the
solution bases
$I^+_X|_{S_{+}}$, $I^-_X|_{S_{-}}$, $e^{2\pi i\alpha_X}I^+_X|_{S_+}$ respectively
form a normalized tuple of sectorial solution bases (see Definition \ref{defnbas}).
\end{proposition}

The Proposition is proved in 2.3.

At the end of the paper we also prove the following more precise asymptotic
statement on the solutions $I_X^{\pm}$. We will not use it in the paper.

\begin{proposition}\label{asymptotic}
For any vertex $X$, the function $I^{\pm}_X(\lambda)$ is a solution of
(\ref{ODE}), with the asymptotic behavior (uniform in the sector $S_{\pm}$)
$$
I^{\pm}_X (\lambda)\sim_{\vert \lambda\vert \to \infty}
D_{X,X}\left( \prod_{j\in X} \Gamma(\alpha_j)\right) e^{-\lambda
f_0(X)} \lambda^{-\alpha_X} v_X,
$$
where $(v_X)_{X\in \xxx}$ is the standard base of $\R^\xxx$,
$$
D_{X,X}=\left(\prod_{j\in X} \Gamma(\alpha_j)\right) \prod_{j\in
X} \vert l_j(e_{X\setminus \{j\}})\vert^{-\alpha_j}
\prod_{j\not\in X} \vert f_j(X)\vert^{\alpha_j}.
$$
\end{proposition}

\renewcommand{\H}{\mathcal H}
\def\rr{\mathbb R}
\subsection{The plan of the computaton of Stokes operators}
For the proof of Theorem  \ref{tstokes} we have to calculate
 the transition matrices $C_0$, $C_1$ between the
 sectorial solution bases from Proposition \ref{canbas}. One has
 \begin{equation}(I^-_X)(\lambda)=(I^+_X)(\lambda)C_0 \ \text{for} \
 \lambda\in\mathbb R_+.\label{stokco}\end{equation}
 This follows from definition and the last statement of Proposition \ref{canbas}.

 To calculate $C_0$, the strategy is
 to pass through the integrals $I_\Delta(\lambda)$, $\Delta\in\mathcal D^{+}$,
 which are well-defined on the axis $\lambda\in\mathbb R_+$.

 For any $\Delta,\Delta'\in\ddd^+$ denote
\begin{equation}
\H(\Delta,\Delta')=\{\text{the hyperplanes} \ H_j \ \text{separating} \
\Delta \ \text{from} \ \Delta'\},\label{separ}
\end{equation}
$$ |\H(\Delta,\Delta')|= \ \text{the cardinality of} \ \H(\Delta,\Delta').$$

\begin{lemma}\label{decomposition}
For $\Re(\lambda)>0$, we have
\begin{equation}
I^+_X(\lambda)=\sum_{\Delta\subset \ccc^+_X}\eta(X,\Delta)
I_\Delta (\lambda), \ \text{and}
\label{i+d}\end{equation}
\begin{equation}
I^-_X(\lambda)=\sum_{\Delta\subset \ccc^+_X}
\overline\eta(X,\Delta) I_\Delta (\lambda),
\label{i-d}\end{equation}
where
$$\eta(X,\Delta)=1, \ \text{if} \ \Delta=\Delta_X, \ \text{otherwise,} \
\eta(X,\Delta)=e^{\pi i\alpha_{\H(\Delta,\Delta_X)}}.$$
\end{lemma}

The Lemma is proved below.

To calculate $C_0$, we have to express $(I^-_X)$ via $(I^+_X)$. The previous
Lemma expresses $I^{\pm}_X$ via the integrals $I_{\Delta}$. Lemma \ref{inversion}
formulated in Section  3
provides the inverse expression of the integrals $I_{\Delta}$ via $I^+_X$.
Afterwards $C_0$ is calculated by substituting the latter inverse expression to (\ref{i-d}).

The proof of Lemma \ref{inversion} is based on the next purely combinatorial
identity, which holds for arbitrary generic arrangement of hyperplanes and a linear
function. To state it, let us introduce some more notations.

For any domain $\Delta\in\ddd^+$ denote
$$\hat\Delta=\Delta\cup(\partial \Delta\cap\partial \mathcal C^+_{X(\Delta)}). \
\text{By definition, for any vertex} \ X\in\xxx$$
\begin{equation}\text{the closed cone} \ \overline{\ccc_X^+} \
\text{is the disjoint union of the sets} \
\hat\Delta, \ \Delta\in\ddd^+, \ \Delta\subset\ccc_X^+.\label{union}\end{equation}

Recall that for any subset $B\subset\rr^n$ $\chi_B:\rr^n\to\rr$ denotes the characteristic
function of $B$: $\chi_B(x)\equiv1$ on $B$, $\chi_B(x)\equiv0$ on $\rr^n\setminus B$.

\begin{lemma} \label{lcomb} Consider arbitrary generic hyperplane arrangement and
a linear function, as at the beginning of the paper. Let $\Delta_X$, $\ccc^+_X$ be the
corresponding domains and cones defined in the Introduction.
For any vertex $X$ of the arrangement one has
\begin{equation}\chi_{\hat \Delta_X}=\sum_{X'\in\partial \Delta_X}\nu(\Delta_X,X')
\chi_{\overline{C_{X'}^+}},
\ \text{where} \ \nu(\Delta,X')=(-1)^{|\H(\Delta,\Delta_{X'})|},\label{chi}\end{equation}
$\H(\Delta,\Delta_{X'})$ was defined in (\ref{separ}).
 \end{lemma}

This Lemma is proved in Section 3.

A version of Lemma \ref{lcomb} was stated and proved by A.N.Varchenko
and I.M.Gelfand in
\cite{Varch3}. Namely they had  shown that the characteristic function of a domain
$\Delta_X$ can be uniquely presented as a  linear combination (with coefficients
$\pm1$) of characteristic functions of some cones (of maybe different dimensions).
They provided some implicit description of  the coefficients of this linear
combination without an explicit formula.
Lemma \ref{lcomb} provides an explicit  formula. Its proof uses a method different
from that of \cite{Varch3}.

\begin{proof} {\bf of Lemma \ref{decomposition}.}
Let us prove formula (\ref{i+d}) of the Lemma. Formula (\ref{i-d}) then follows
from (\ref{i+d}), the equality
\begin{equation}I^-_X(\lambda)=\overline{I^+_X(\lambda)} \ \text{for any} \
\lambda\in\mathbb R_+\label{ipmbar}\end{equation}
and the complex conjugatedness of the right-hand sides of (\ref{i+d}) and (\ref{i-d})
(the integrals $I_{\Delta}(\lambda)$ are real for $\lambda\in\mathbb R_+$).
It suffices to show that
the analytic extension of the integrals $I^+_X$ to the semiaxis $\lambda>0$
is defined by formula (\ref{i+d}). Indeed, the mapping
$F_{\rho}:x\mapsto X+\rho(x-X)$
is a real-linear isomorphism $\ccc^+_X\to\ccc^{\rho}_X$ that tends to the
identity, as $\rho\to1$. The cone
$\ccc^{\rho}_X$ is the union of the closures of the domains
$\Delta(\rho)=F_{\rho}(\Delta)$, $\Delta\in\ddd^+$,
$\Delta\subset\ccc^+_X$. The integral $I^+_X$ is the sum of the integrals
of $e^{-\lambda f_0}\Omega|_{\ccc^{\rho}_X}$ over the domains $\Delta(\rho)$.
Each latter integral tends (as $\rho\to1$) to the integral over
$\Delta$ of the form $e^{-\lambda f_0}\Omega$, where the branch $\Omega|_{\Delta}$
is the immediate analytic extension
 of $\Omega|_{\cup_{\Im\rho<0}\Delta(\rho)}$ to $\Delta=\Delta(1)$.
We claim that thus extended branch $\Omega|_{\Delta}$ is the
$\eta(X,\Delta)$- th multiple of the standard real branch of $\Omega$
on $\Delta$ (see (\ref{int})).
Indeed, fix a $x\in\Delta$ and denote $L\subset\C^N$ the
complex line containing the segment $[X,x]$.
The latter segment
intersects $\Delta_X$ (by definition; fix a point of their intersection and denote
it $x_0$). Fix a $\rho\in\mathbb C$, $|\rho|=1$, with $\Im\rho<0$. Denote
$$x_{\rho}=F_{\rho}(x)\in\Delta(\rho), \ \delta_{\rho}=
\{ F_{e^{i\theta}}(x) \ | \arg\rho\leq\theta
\leq0\}. \ \text{Consider the path}$$
$$\gamma:[0,1]\to\C^N \ \text{from} \ x_0 \ \text{to} \ x: \
\gamma=[x_0,x_{\rho}]\circ\delta_{\rho}.$$
By construction, the previously constructed branch $\Omega|_{\Delta}$
is obtained by the analytic extension of the standard real branch of
$\Omega$ on $\Delta_X$ along the path $\gamma$ (all the points of $\gamma$ except
for its ends $x_0$ and $x$ lie in $\cup_{\im\rho<0}\ccc^{\rho}_X$).
For any hyperplane $H_j$ intersecting the segment $[x_0,x]$ denote
$x_j$ the intersection point.
We consider that the point $x$ is chosen generic so that the points $x_j$
are distinct. The path $\gamma$ is
isotopic in $L\setminus\cup_jH_j$ to the segment $[x_0,x]$ where small
intervals $(a_j,b_j)$ containing $x_j$ are replaced by half-circles in $L$
(with the same ends $a_j$ and $b_j$) oriented
counterclockwise (the notion "counterclockwise" is independent on the choice of affine
complex coordinate on $L$). Extending the form $\Omega$ along a previous
half-circle yields extra multiplier $e^{\pi i\alpha_j}$.
This implies that the extended branch $\Omega|_{\Delta}$ is
the standard real branch times $\eta(X,\Delta)$. This together with the
previous discussion proves the Lemma.
\end{proof}

\subsection{The integrals $I_{X}^{\pm}$. Proof  of Proposition \ref{canbas}.}
The vector functions
$I^{\pm}_X$ are linear combinations of integrals over domains $\Delta$
(Lemma \ref{decomposition}). Therefore, they are solutions of (\ref{ODE}),
as are the latter integrals (see the Introduction). Now we have
to show that they form canonical sectorial solution bases.

Given a ray in $\C$, we say
that a vector function
$f$ is {\it asymptotically bigger} than another one $g$  along the ray, if
 $g(z)=o(f(z))$, as $z\to\infty$ along the ray. A collection of functions is
 {\it asymptotically ordered} along a ray, if for any two distinct
  functions one is asymptotically  bigger than the other one.
 We use the following characterization of canonical solution bases, which
 follows from the general theory of linear equations with irregular
 singularities.
 \begin{proposition} Consider arbitrary canonical solution base
 of (\ref{ODE}) in $S_+$ (or in $S_-$). The basic solutions are asymptotically
 ordered along both semiaxes $\pm\lambda>0$; their orderings along
 these semiaxes are opposite to each other.
 Vice versa, given arbitrary collection of solutions $F_X$ (numerated by
 all the vertices $X$)
 of (\ref{ODE}) in the sector $S_{\pm}$ under consideration. Let $F_X$ be asymptotically
 ordered along the previous semiaxes and the corresponding orderings be opposite to
 each other.  Then $F_X$ is a canonical sectorial solution base.
 \end{proposition}

 {\bf Addendum.} {\it Let $S_{\pm}$, $S_0$, $S_1$, $S_2$ be as in (\ref{sectors}) and
 (\ref{nsectors}).
 Let $F^{\pm}=(F^{\pm}_{X})_{X\in\xxx}$ be a pair of canonical sectorial solution bases  in $S_{\pm}$ such that for any $X\in\xxx$ one has
 \begin{equation}F^-_{X}(\lambda)=F^+_{X}(\lambda)+o(F^+_{X}(\lambda)), \
 \text{as} \ \lambda\in\mathbb R_+,
  \ \lambda\to+\infty.\label{1+o}\end{equation}
 Then the liftings to $S_0$, $S_1$, $S_2$ of the  bases
 $F^+|_{S_+}$, $F^-|_{S_-}$, $(e^{2\pi i\alpha_X}F^+_X)|_{S_+}$ form
 a normalized tuple of canonical sectorial bases.}

 \medskip

 \begin{proof} The statements of the Proposition and the Addendum are obvious for the
 formal normal form  (\ref{nform}). Let us prove the statements of the Addendum for
 (\ref{nform}) in more detail. Each solution base of (\ref{nform}) under consideration is
 defined by a diagonal fundamental matrix. Any two (locally defined) diagonal
 fundamental matrices are obtained one from the other by multiplication of the
 diagonal elements by appropriate constants. The latter constants comparing the
 fundamental matrices of $F^+|_{S_0}$ and $F^-|_{S_1}$ on $S_{01}$
 ($F^-|_{S_1}$ and $(e^{2\pi i\alpha_X}F^+_X)|_{S_2}$ on $S_{12}$) are unit, i.e.,
 the three latter solution bases are holomorphic extensions of each other. This follows
 from (\ref{1+o}) (for the former base pair) and the fact that the solution base
 $(e^{2\pi i\alpha_X}F^+_X)|_{S_+}$  is the image of $F^+|_{S_+}$ under the
 clockwise
 monodromy around 0. Hence, the lifted bases from the Addendum form  a normalized
 tuple (see Definition \ref{defnbas}).

 Now given arbitrary differential equation (\ref{ODE}). Consider the
 variable transformations
 inverse to its sectorial normalizations. These transformations send (\ref{nform})
 to (\ref{ODE}) and thus, the canonical sectorial solution bases of (\ref{nform}) to those
 of (\ref{ODE}), and preserve the asymptotic orderings and relations (\ref{1+o}).
 This together with the statements of the Proposition and the Addendum for (\ref{nform})
 proves them for (\ref{ODE}).
 \end{proof}

One has
\begin{equation}I^{\pm}_X=I_{\Delta_X}+o(I_{\Delta_X}), \ \text{as} \ \lambda\in
\mathbb R_+, \
\lambda\to+\infty,\label{asint}\end{equation}
This follows from Lemma \ref{decomposition},
the inclusion $\Delta_X\subset\ccc^+_X$
and the inequality $f_0|_{\overline\Delta_{X'}}>f_0(X)$ valid for any vertex
$X'\in\overline\ccc^+_X\setminus X$ (which holds by definition).
The integrals $I_X^{\pm}$ are asymptotically ordered along the semiaxis
$\lambda>0$: $I_X^{\pm}$ is asymptotically greater than $I_{X'}^{\pm}$, if and only if
$f_0(X)<f_0(X')$. This follows from (\ref{asint}) and the previous inequality.
The same integrals $I_X^{\pm}$ are also asymptotically ordered along the opposite
semiaxis $\lambda<0$, and their latter order is opposite to the previous one. Indeed,
let us prove the latter statement for $I^+_X$. Then for $I^-_X$ the same statement
follows from the one for $I^+_X$ and
the relation $I^-_X(\lambda)=\overline{I^+_X(\bar\lambda)}$ (which
follows from (\ref{ipmbar})). The cone $\ccc^{-1}_X=\ccc^-_X$ is adjacent to the union
$\cup_{\im\rho<0}\ccc^{\rho}_X$ and is a union of domains from $\ddd^+$
(denote $\Delta^-_X\subset\ccc^-_X$ the domain with vertex at $X$). The integral
$I^+_X(\lambda)$ restricted to $\lambda\in\mathbb R_-$ can be expressed as a linear
combination of the integrals over the previous domains, as in Lemma
\ref{decomposition} and its proof. The  integral $I^+_{\Delta_X^-}$ appears there with
the coefficient $e^{-\pi i\alpha_X}$. One has
$$I_X^{\pm}(\lambda)=e^{-\pi i\alpha_X}I_{\Delta_X^-}(\lambda)+o(I_{\Delta_X^-}
(\lambda)), \ \text{as} \ \lambda\in\mathbb R_-, \ \lambda\to-\infty,$$
as in (\ref{asint}). This together with the arguments following (\ref{asint}) prove the
previous asymptotic order statement.

The two  asymptotic order statements proved above together with the previous
Proposition imply that the integrals $I^{\pm}_X$ form canonical solution bases in
$S_{\pm}$. This proves the first part of Proposition \ref{canbas}.

Let us prove the second part of  Proposition \ref{canbas} (about the normalized base
tuple).
By the Addendum, to do this, it suffices to prove equality (\ref{1+o}) for the bases
$F^{\pm}_X=I^{\pm}_X$. This equality follows immediately
from (\ref{asint}).  Proposition \ref{canbas} is proved.

\section{The relations between $I_X^\pm$ and $I_\Delta$. Proof of Theorem
\ref{tstokes}}

As it is shown (at the end of the Section),
Theorem \ref{tstokes} is implied by Lemma \ref{decomposition}
and the following Lemma. The proof of the latter is based on Lemma \ref{lcomb};
both Lemmas are proved below.

\begin{lemma} \label{inversion}
For any $\Delta\in\ddd^+$ the following equalities hold for all $\lambda\in\mathbb R_+$:
\begin{equation}
I_\Delta (\lambda)= \sum_{X\in \partial\Delta} \psi(\Delta, X)
I_X^+(\lambda),
\label{inv+}\end{equation}
\begin{equation}
I_\Delta (\lambda)= \sum_{X\in \partial\Delta}
\overline{\psi(\Delta, X)} I_X^-(\lambda), \ \text{where}
\label{inv-}\end{equation}
$$
\psi(\Delta, X)=1 \ \text{if} \ \Delta=\Delta_X, \ \text{otherwise,} \
\psi(\Delta,X)=(-1)^{|\H(\Delta,\Delta_X)|}e^{i\pi\alpha_{\H(\Delta,\Delta_X)}},
$$
the set $\H(\Delta,\Delta_X)$ was defined in (\ref{separ}).
\end{lemma}

\begin{proof} {\bf of Lemma \ref{lcomb}.} Fix a vertex $X$ and denote
 $$\mathcal D^+_X=\{ \Delta\in\mathcal D^+ \ | \ X(\Delta)\geq X\}.$$
 (Recall that the vertices are ordered so that the function
 $X\mapsto f_0(X)$ is increasing.)
 The domain collection $\mathcal D^+_X$ is in 1-to-1 correspondence with
 the vertices $X'\geq X$. (Denote $M$ the number of elements
 in each collection.) Each domain $\Delta\subset\ccc^+_X$ is contained in
 $\mathcal D^+_X$ by definition and since $f_0|_{\ccc^+_X}$ is bounded from below.
 By (\ref{union}), for any vertex  $X'\geq X$ one has
 $$\chi_{\overline{\ccc^+_{X'}}}=\sum_{\Delta\subset\ccc^+_{X'}}\chi_{\hat\Delta}=
 \sum_{\Delta\in\mathcal D^+_X}\theta(X',\Delta)\chi_{\hat\Delta}, \ \text{where}$$
 $$\theta(X',\Delta)=1 \ \text{whenever} \ \Delta\subset\ccc^+_{X'}; \ \theta(X',\Delta)=0 \ \text{otherwise}.$$
 In other terms, the vector of the functions $\chi_{\overline{\ccc^+_{X'}}}$ is obtained
 from the vector of the functions $\chi_{\hat\Delta}$ by multiplication by the
 $M\times M$ matrix $\theta(X',\Delta)$ with indices
 $X'\geq X$  and $\Delta\in\mathcal D^+_X$.

 For the proof of (\ref{chi}) we extend the
 values $\nu(\Delta,X')$ (which were defined in (\ref{chi}) for $X'\in\partial\Delta$)
 up to a $M\times M$- matrix (with the previous indices) by putting
 $$\nu(\Delta,X')=0 \ \text{whenever} \ X'\notin\partial\Delta.$$
We show that the matrices $\nu(\Delta,X')$ and $\theta(X',\Delta)$ are inverse,
 i.e., for any two vertices $X',X''\geq X$ one has
 \begin{equation}\sum_{\Delta\in\mathcal D^+_X}\theta(X',\Delta)\nu(\Delta,X'') \
 \text{equals} \ 0 \ \text{if} \
 X'\neq X'' \ \text{and equals} \ 1 \ \text{if} \ X'=X''.\label{obr}\end{equation}
 This will prove the Lemma.

 The only nonzero terms of the sum in (\ref{obr}) correspond exactly to
 $\Delta\in D(X',X'')$, where
 \begin{equation}D(X',X'')=\{\Delta\subset\ccc^+_{X'} \ | \ X''\in\partial\Delta\}; \
 \text{one has} \  X'\leq X'', \ \text{if} \ D(X',X'')\neq\emptyset.
 \label{incl}\end{equation}

 Case $X'=X''$.  \ Then \ $D(X',X'')=\{\Delta_{X'}\}$ \ \ and  \ \
  $\theta(X',\Delta_{X'})=\nu(\Delta_{X'},X')=1$ by definition. This proves the
  second statement of (\ref{obr}).

 Case $X'>X''$. Then all the terms of the sum in (\ref{obr}) vanish, see
 (\ref{incl}).

  Case $X'<X''$.
  Let us introduce affine coordinates $x_1,\dots,x_n$
  on $\rr^n$ so that $X''$ is the origin and
  the arrangement hyperplanes through $X''$ are the coordinate hyperplanes.
  Fix a hyperplane $H=\{ x_j=0\}$  (which contains $X''$)
  that does not contain  $X'$ (it exists by definition).

   If $X''\in\ccc^+_{X'}$, then the domains $\Delta\in D(X',X'')$ intersect a small
   neighborhood of $X''$ by  the coordinate quadrants
  (whose number equals $2^n$). If $X''\in\partial\ccc^+_{X'}$, then locally near $X''$
  the cone $\ccc^+_{X'}$ is the coordinate cone defined by the inequalities
  $\pm x_j>0$ (for a certain collection of distinct indices $j\neq i$); the  domains
  $\Delta\in D(X',X'')$ are locally the coordinate quadrants in the latter cone.
  In both cases the domain collection $D(X',X'')$  is split into pairs. The
 domains in each pair are {\it adjacent across} $H$: by definition, this means
  that  they are adjacent to a common face in $H$ (of the same dimension, as $H$),
  and thus, are separated from each other by $H$. For any two domains
  $\Delta_1$ and $\Delta_2$ adjacent across $H$ one has
  $\nu(\Delta_1,X'')+\nu(\Delta_2,X'')=0$ (hence, the corresponding terms of the
  sum in (\ref{obr}) cancel out and the latter sum vanishes). Indeed, let
  $H$ separate $\Delta_1$ from $\Delta_2$ and $\Delta_{X''}$ (otherwise we interchange
  $\Delta_1$ and $\Delta_2$). Then
  $$\H(\Delta_1,\Delta_{X''})=\H(\Delta_2,\Delta_{X''})\cup H$$
  by definition. This together with the definition of $\nu(\Delta_j,X'')$, see (\ref{chi}),
   proves the previous cancellation statement, (\ref{obr}) and Lemma \ref{lcomb}.
   \end{proof}

\begin{proof} {\bf of Lemma \ref{inversion}.} Let us prove (\ref{inv+}) (then
(\ref{inv-}) follows by complex conjugation argument, see (\ref{ipmbar})).
Let us substitute the expression (\ref{i+d}) for $I^+_X$ via the integrals over domains
to the right-hand side of (\ref{inv+}). We show
that for any $\Delta'\in\ddd^+$ the corresponding coefficients at $I_{\Delta'}$
obtained by this substitution cancel out, except for the unit coefficient
corresponding to $\Delta'=\Delta$. This will prove the Lemma. After the
previous substitution the right-hand side of (\ref{inv+}) takes the form
$$\sum_{X\in\partial\Delta}\sum_{\Delta'\subset\ccc^+_X}
\eta(X,\Delta')\psi(\Delta,X)I_{\Delta'}, \ \eta(X,\Delta') \ \text{are the same, as in
(\ref{i+d})}.$$
For any $X$, $\Delta$, $\Delta'$ such that $X\in\partial\Delta$, $\Delta'\subset\ccc^+_X$
one has
\begin{equation}\eta(X,\Delta')\psi(\Delta,X)=(-1)^{|\H(\Delta_X,\Delta)|}
e^{i\pi\alpha_{\H(\Delta,\Delta')}}.\label{heta}\end{equation}
Indeed, recall that by definition,
\begin{equation}\eta(X,\Delta')=e^{\pi i\alpha_{\H(\Delta',\Delta_X)}}, \
\psi(\Delta,X)=(-1)^{|\H(\Delta_X,\Delta)|}e^{i\pi \alpha_{\H(\Delta_X,\Delta)}}.
\label{hetanew}\end{equation}
Formula (\ref{heta}) follows from (\ref{hetanew})  and the fact that for any
$\Delta\in\ddd^+$, $X\in\partial\Delta$ and $\Delta'\subset\ccc^+_X$ one has
\begin{equation}\H(\Delta,\Delta_X)\cap \H(\Delta_X,\Delta')=\emptyset,
\ \H(\Delta,\Delta_X)\cup \H(\Delta_X,\Delta')=\H(\Delta,\Delta').\label{hincl}
\end{equation}
Indeed, each hyperplane $H\in\H(\Delta_X,\Delta')$, which separates $\Delta_X$ from
$\Delta'$, by definition, also separates $\Delta$ from $\Delta'$. Otherwise
$H$ separates $\Delta$ from $\Delta_X$ (hence, $X\in H$). Therefore, $H$
does not cut the cone $\ccc_X^+$ and thus, cannot separate its subdomains
$\Delta_X$ and $\Delta'$, - a
contradiction. Each $H\in \H(\Delta,\Delta_X)$ separates $\Delta$ from $\Delta'$, since
it separates $\Delta$ from the cone $\ccc^+_X\supset\Delta'$ (which
follows from definition). Thus,
$$\H(\Delta,\Delta_X)\cup\H(\Delta_X,\Delta')\subset
\H(\Delta,\Delta').$$
Vice versa, each hyperplane $H\in \H(\Delta,\Delta')$ separates
$\Delta$ from $\Delta'$ (by definition), and $\Delta_X$ is either on the $\Delta'$- s
or on the $\Delta$- s side. These two (incompatible) cases take place, when
$H\in\H(\Delta,\Delta_X)$ (respectively, $H\in\H(\Delta_X,\Delta')$). This proves
(\ref{hincl}) and (\ref{heta}).

Now by (\ref{heta}), the right-hand side of (\ref{inv+}) equals the linear combination of
the integrals $I_{\Delta'}$ with the coefficients
$$e^{i\pi\alpha_{\H(\Delta,\Delta')}}\sum_{X\in\partial\Delta, \
\Delta'\subset\ccc_X^+}
(-1)^{|\H(\Delta_X,\Delta)|}.$$
The latter sum over vertices $X$ equals the value on $\Delta'$ of the characteristic
function combination
(\ref{chi}) (with $\Delta_X$, $X'$ in (\ref{chi}) replaced by $\Delta$, $X$ respectively) by definition. Hence, it vanishes, if $\Delta'\neq\Delta$, and equals 1
if $\Delta'=\Delta$ (Lemma \ref{lcomb}). This proves Lemma \ref{inversion}.
\end{proof}

\begin{proof} {\bf of Theorem \ref{tstokes}.} Let $C_0=(C_0(X',X))_{X',X\in\xxx}$ be the
Stokes matrix (\ref{stokes}) corresponding to the normalized base tuple in $S_0$,
$S_1$, $S_2$ from Proposition \ref{canbas}. One has
\begin{equation}I^-_X(\lambda)=\sum_{X'\in\xxx}C_0(X',X)I^+_{X'}(\lambda)
 \ \text{for all} \ \lambda\in\mathbb R_+,\label{expst}\end{equation}
by definition. Let us calculate the coefficients $C_0(X',X)$.
Lemma \ref{decomposition} gives formula (\ref{i-d}) for $I^{-}_X$ as a linear
combination of the integrals $I_{\Delta}$ with constant coefficients.
Replacing each $I_{\Delta}$ in (\ref{i-d}) by its expression (\ref{inv+}) via the integrals
$I^+_{X'}$ yields (\ref{expst}) with
\begin{equation}
C_0(X',X)=\sum_{\Delta\subset\ccc^+_X, \ X'\in\partial\Delta}\gamma(X,X',\Delta),
\label{cosum}\end{equation}
\begin{equation}
\gamma(X,X',\Delta)=\overline\eta(X,\Delta)\psi(\Delta,X')=
(-1)^{|\H(\Delta,\Delta_{X'})|}
e^{\pi i(\alpha_{\H(\Delta,\Delta_{X'})}-\alpha_{\H(\Delta,\Delta_X)})}.\label{stokg}
\end{equation}
In the case, when $X'=X$, obviously $C_0(X',X)=1$. If $X'\notin\overline\ccc^+_X$, then
$C_0(X',X)=0$, since the previous sum contains no terms.

Thus, everywhere
below in the calculation of $C_0$ we consider that  $X'\in\overline\ccc^+_X\setminus X$.
Let us calculate the  sum (\ref{cosum}). To do this,
we extend (literally) the definition
of $\H(\Delta_1,\Delta_2)$ to the case, when each
$\Delta_j$ is an arbitrary union of domains in $\ddd^+$, by putting
$\H(\Delta_1,\Delta_2)$  to be the number of the arrangement
hyperplanes separating $\Delta_1$ from $\Delta_2$. Then we extend analogously
the definition of the values $\gamma(X,X',\Delta)$ (for $\Delta$ being a union of
domains) by writing formula (\ref{stokg}) with thus generalized $\H(\Delta,\Delta_{X'})$,
$\H(\Delta,\Delta_X)$.

Fix an arbitrary arrangement hyperplane $H_j$ through $X'$ that does not contain $X$
and a pair of domains $\Delta_1,\Delta_2\subset\ccc^+_X$ adjacent across $H_j$
(see the proof of Lemma \ref{lcomb} in the previous Subection), $X'\in\partial\Delta_l$,
$l=1,2$. Let us compare the values $\gamma(X,X',\Delta_l)$.

Case 1: the pair $(X,X')$ is positive exceptional
and the hyperplane $H_j$ is exceptional
(see Definition \ref{dexcep}; then $X'\in\mathcal C_X^+$ and no arrangement hyperplane
through $X'$ contains $X$; thus, $H_j$ can be chosen arbitrary, e.g., exceptional). We
claim that
\begin{equation}\gamma(X,X',\Delta_1)+\gamma(X,X',\Delta_2)=0.
\label{ggo}\end{equation}
Indeed, by definition, the domains $\Delta_X$ and $\Delta_{X'}$ lie on the same side
from $H_j$. Let $\Delta_1$ also lie on the same side; then $\Delta_2$ lies on the
other side (otherwise, we interchange $\Delta_1$ and $\Delta_2$). One has
$$\H(\Delta_2,\Delta_X)=\H(\Delta_1,\Delta_X)\cup H_j,
\ \H(\Delta_2,\Delta_{X'})=\H(\Delta_1,\Delta_{X'})\cup H_j,$$
since $H_j$ is the only arrangement hyperplane separating $\Delta_1$ and $\Delta_2$.
This together with (\ref{stokg}) implies (\ref{ggo}).

Case 2: the pair $(X,X')$ is not positive exceptional.
(This includes the case, when $X'\in\partial\ccc^+_X$, since then any hyperplane
through $X'$ that does not contain $X$ (thus, $H_j$) separates $\Delta_{X'}$ from
$\Delta_X$. This follows from definition and the increasing of the   function $f_0$ along
the segment $[X,X']$ oriented from $X$ to $X'$.)  We claim that
\begin{equation}
\gamma(X,X',\Delta_1)+\gamma(X,X',\Delta_2)=-(2i\sin\pi\alpha_j)
\gamma(X,X',\Delta_1\cup\Delta_2),\label{uniond}\end{equation}
and this equality remains valid in the case, when $\Delta_1$ and $\Delta_2$ are
{\it adjacent across $H_j$ unions
of domains} in $\ccc^+_X$. The latter means that the
domains from $\Delta_1$, $\Delta_2$  have the following properties:

1)  the closure of each domain in $\Delta_1$, $\Delta_2$ contains $X'$;

2) each domain in $\Delta_1$ is adjacent across $H_j$
to a domain in $\Delta_2$ and vice versa.

Indeed, without
loss of generality we consider that $\Delta_1$, $\Delta_{X'}$
are separated by $H_j$ from $\Delta_2$ and $\Delta_X$ (interchanging
$\Delta_1$ and $\Delta_2$ if necessary).  By definition, one has
$$\H(\Delta_2,\Delta_{X'})=\H(\Delta_1,\Delta_{X'})\cup H_j, \
\H(\Delta_1,\Delta_{X'})=\H(\Delta_1\cup\Delta_2,\Delta_{X'}),$$
$$\H(\Delta_1,\Delta_X)=\H(\Delta_2,\Delta_X)\cup H_j, \
\H(\Delta_2,\Delta_X)=\H(\Delta_1\cup\Delta_2,\Delta_X).$$
Hence, by (\ref{stokg}),
$$\gamma(X,X',\Delta_1)=e^{-\pi i\alpha_j}
 \gamma(X,X',\Delta_1\cup\Delta_2),$$
$$ \gamma(X,X',\Delta_2)=-e^{\pi i\alpha_j}\gamma(X,X',\Delta_1\cup\Delta_2).$$
The two latter formulas imply (\ref{uniond}).

If the pair $(X,X')$ is positive exceptional, then $C_0(X',X)=0$. Indeed, fix an
exceptional hyperplane $H_j$. The collection of all the domains in $\ccc^+_X$ whose
closures contain $X'$ is split into pairs of adjacent domains across $H_j$.
The terms in the sum (\ref{cosum}) corresponding
to two adjacent domains cancel out by (\ref{ggo}), hence the sum vanishes.

Let now the pair $(X,X')$ be not positive exceptional. Let  us numerate all the
hyperplanes $H_{j_1},\dots,H_{j_q}$ through $X'$ that do not contain $X$
(one has $q\leq k$). If $X'\in\ccc^+_X$, then $q=k$ and these are all the arrangement
hyperplanes through $X'$.
Otherwise, if $X'\in\partial \ccc^+_X$, then $q<k$ and these are all the arrangement
hyperplanes through $X'$ that do not contain $X$ (or equivalently, that do not contain
faces of the cone $\ccc^+_X$). In both cases one has $\{ j_1,\dots,j_q\}=X'\setminus X$.
The  terms in the sum (\ref{cosum})  correspond
to the domains $\Delta_1,\dots\Delta_{2^q}$, which we numerate as follows.
Put $\Delta_1=\Delta_{X'}$, $\Delta_2$ be the domain adjacent across $H_{j_1}$ to
$\Delta_1$, $\Delta_3$ ($\Delta_4$) be the domain adjacent across $H_{j_2}$ to
$\Delta_1$ (respectively, $\Delta_2$), etc., for any $s=1,\dots,q-1$ the domains
$\Delta_{2^{s}+1},\dots\Delta_{2^{s+1}}$ are adjacent across $H_{j_{s+1}}$ to
$\Delta_1,\dots,\Delta_{2^s}$. We claim that for any $s=1,\dots,q$
\begin{equation} \sum_{l=1}^{2^s}\gamma(X,X',\Delta_l)=
\gamma(X,X',\cup_{l=1}^{2^s}\Delta_l)\prod_{r=1}^s(-2i\sin\pi\alpha_{j_r}),
\label{union1}\end{equation}
\begin{equation} \sum_{l=2^s+1}^{2^{s+1}}\gamma(X,X',\Delta_l)=
\gamma(X,X',\cup_{l=2^s+1}^{2^{s+1}}\Delta_l)\prod_{r=1}^s(-2i\sin\pi\alpha_{j_r}),
\ \text{whenever} \ s<q.
\label{union2}\end{equation}
We prove both statements (\ref{union1}), (\ref{union2}) by induction in $s$.

The induction base for $s=1$ follows from (\ref{uniond}) and the fact that
$\Delta_3$, $\Delta_4$ are adjacent across $H_{j_1}$ (by definition).

Induction step. Let (\ref{union1}), (\ref{union2}) be proved for a given $s<q$.
Let us prove (\ref{union1}) for $s$ replaced by $s+1$. The domain unions from
(\ref{union1}) and (\ref{union2}) are adjacent across $H_{j_{s+1}}$ to each other by
definition. Adding  equalities (\ref{union1}) and (\ref{union2}) and applying
(\ref{uniond}) to the $\gamma$'s in the right-hand side yields (\ref{union1})
for $s$ replaced by $s+1$. Equality (\ref{union2}) for $s+1\leq q$ is proved analogously.
The induction step is over and statements (\ref{union1}), (\ref{union2}) are proved.

 Formula (\ref{union1}) with $s=q$ says that the sum (\ref{cosum}) equals
$$\gamma(X,X',\wt\Delta)\prod_{s=1}^q(-2i\sin\pi\alpha_{j_s}), \ \text{where}$$
$\wt\Delta$ is the union of all the domains in $\ccc^+_X$ whose closures contain $X'$.
The latter expression coincides with
the right-hand side in (\ref{co}), by (\ref{stokg}) (applied to $\widetilde\Delta$) and since
$$A=\H(\wt\Delta,\Delta_X), \ B=\H(\wt\Delta,\Delta_{X'}), \ q=|X'\setminus X|$$
(by definition). This proves (\ref{co}).

\def\odlo{\text{(\ref{ODE})}(f_0)}
\def\odmlo{\text{(\ref{ODE})}(-f_0)}

Now let us prove (\ref{c1}). The Stokes matrix $C_1$ is the transition
matrix between the canonical solution bases $I^-_X(\lambda)$ and
$e^{2\pi i\alpha_X}I^+_X(\lambda)$, $\lambda\in\mathbb R_-$, by definition and
Proposition \ref{canbas}.
To calculate it, we consider the variable change $\lambda\mapsto-\lambda$,
which transforms the equation $\text{(\ref{ODE})}=\text{(\ref{ODE})}(f_0)$ to
the new one (denoted $\text{(\ref{ODE})}(-f_0)$). The latter equation
corresponds to  the same hyperplane arrangement equipped with the new
linear function
$$\wt f_0=-f_0.$$
Denote $J^{\pm}_X(\lambda)$ the canonical basic solutions of $\odmlo$ in the sector
$S_{\pm}$: the solutions given by Proposition \ref{canbas} (denoted there by
$I^{\pm}_X(\lambda)$).
The variable change $\lambda\mapsto-\lambda$ transforms the canonical sectorial
basic solutions of $\odlo$  in $S_{\pm}$ to those of $\odmlo$ in $S_{\mp}$. We show that
\begin{equation}I^-_X(-\lambda)=e^{\pi i\alpha_X}J^+_X(\lambda) \ \text{for all} \
\lambda\in S_+.\label{ijl}\end{equation}
Then one has
\begin{equation} e^{\pi i\alpha_X}J^-_X(\lambda)=\sum_{X'}C_1(X',X)e^{\pi i\alpha_{X'}}
J_{X'}^+(\lambda), \ \lambda\in\mathbb R_+. \label{defc1}\end{equation}
This follows from definition, (\ref{ijl}) and formula
$$e^{2\pi i\alpha_X}I^+_X(-\lambda)=e^{\pi i\alpha_X}J^-_X(\lambda) \ \text{for any} \
\lambda\in S_-.$$
The latter formula follows from (\ref{ijl}), the fact that $I^-_X|_{S_1}$,
$e^{2\pi i\alpha_X}I^+|_{S_2}$ form a normalized base pair, as do $J^+_X|_{S_0}$,
$J^-_X|_{S_1}$ (Proposition \ref{canbas} applied to  $\text{(\ref{ODE})}(f_0)$
 and $\text{(\ref{ODE})}(-f_0)$),
and Remark \ref{rnbas}. Formula (\ref{defc1}) together with the (already proved)
formula (\ref{co}) for the transition matrix between $J^+_X$ and
$J^-_X$ yields (\ref{c1}). (Here $B^+$ and "positive exceptional" are replaced by
$B^-$ and "negative exceptional", since the sign of the function $f_0$
(which defines the cone $\ccc^+_X$) is changed.)

Let us prove (\ref{ijl}). Let $\ccc^{\rho}_X$, $\rho\in\mathbb C$, $|\rho|=1$, be the
cones defined in (\ref{crho}). By definition,
\begin{equation}J^+_X=(J^{+}_{X,X'})_{X'\in\xxx}, \ J^{+}_{X,X'}=
J^{\rho}_{X,X'} (\lambda)=\int_{\ccc_X^{\rho}} e^{\lambda f_0(x)}\Omega_{X'},
\label{j+}\end{equation}
 \begin{equation}I^-_{X,X'}(-\lambda)=I^{\rho}_{X,X'}(-\lambda)=
\int_{\ccc_X^{\rho}} e^{\lambda f_0(x)}
\Omega_{X'}; \ \Im\rho>0, \ \Re(\rho\lambda)<0.\label{i+}\end{equation}
In  formulas (\ref{j+}) (respectively, (\ref{i+})) the analytic branch of $\Omega_{X'}$
(denoted $\Omega_{X'}^+$ (respectively, $\Omega_{X'}^-$))
in the union $\hat C=\cup_{\Im\rho>0}\ccc_X^{\rho}$ is defined as a result of immediate
analytic extension
of its standard real branch in a neighborhood of $X$ in $\ccc^-_X=\ccc^{-1}_X$
(respectively, $\ccc^+_X=\ccc^1_X$) to the latter union.   One has
\begin{equation}\Omega_{X'}^-=e^{\pi i\alpha_X}\Omega_{X'}^+.\label{opm}
\end{equation}
(This together with (\ref{j+}) and (\ref{i+}) implies (\ref{ijl}).) Indeed,
consider a point $x_0\in\Delta_X\subset\ccc^+_X$ and a path
$$\Gamma:[0,1]\to \hat C, \ \Gamma(t)=X+e^{i\pi t}(x_0-X);$$
$x_0$ being close enough to $X$ in order that $\Gamma(1)\in\ccc^-_X$ be
not separated from $X$ by arrangement hyperplanes. The result of the analytic
extention of $\Omega^-_{X'}$ from $x_0$ along $\Gamma$ is $e^{i\pi\alpha_X}$
times the real branch of $\Omega_{X'}$ defined near $\Gamma(1)$. The latter
branch equals $\Omega^+_{X'}$ by definition. This proves (\ref{opm}) and
hence (\ref{ijl}). The proof of Theorem \ref{tstokes} is complete.
\end{proof}

\section{Appendix: the differential equation (\ref{ODE})}
The proof of (\ref{ODE}) is based on two types of relation; the
first one comes from the fact that $f_0$ and $(f_j)_{j\in X}$ are
linked for any vertex $X$. Indeed, since $f_0-f_0(X)$ and
$(f_j)_{j\in X}$ vanishes at the point $X$, it implies that there
exits constants $(c_{0,j})_{j\in X}$ such that
$$
f_0(z)=f_0(X)+\sum_{j\in X} c_{0,j} f_j (z), \;\;\; \forall z\in
\R^k.
$$
The second relation is of a cohomological type. Let
$U=\{j_1,\ldots ,j_{k-1}\}$ and
$$
\w_U={df_{j_1}\slash f_{j_1}}\wedge \cdots \wedge
{df_{j_{k-1}}\slash f_{j_{k-1}}},
$$
where the points of $U$ are ordered so that the form
$$
df_0\wedge {df_{j_1}}\wedge \cdots \wedge {df_{j_{k-1}}}
$$
is positively oriented. We have
\begin{eqnarray*}
&& d\left( e^{-\lambda f_0}\left(\prod_{j} \vert
f_j\vert^{\alpha_j}\right) \w_U\right)
\\
&=& \left( e^{-\lambda f_0}\prod_{j} \vert
f_j\vert^{\alpha_j}\right) \left(-\lambda df_0\wedge
\w_U+\sum_{j\in U^c} \alpha_j {df_j\slash f_j}\wedge \w_U\right)
\end{eqnarray*}
We see that the orientation of $df_j\wedge {df_{j_1}}\wedge \cdots
\wedge {df_{j_{k-1}}}$ depends on the relative orientation of the
linear forms $df_j$ and $df_0$ on the edge $L_U$. More precisely,
its orientation is equal to the sign of $l_j(e_U)$ (where $e_U$ is
defined in (\ref{eU}), and $l_j=f_j-f_j(0) $ is the linear form
associated with $f_j$). Hence, we have
$$
{df_j\slash f_j}\wedge \w_U= \epsilon (j,U) \w_{U\cup\{j\}},
$$
where
$$
\epsilon (j,U)=\sgn(l_j(e_U)).
$$
To apply Stockes, we need to prove that the boundary terms do not
contribute. Since the integrant $ e^{-\lambda f_0}\prod_{j} \vert
f_j\vert^{\alpha_j} \w_U$ may diverge on the boundary we first
apply Stockes in the subdomain $\Delta^\eta$ defined as follows:
let $\epsilon_i^\Delta $ be the sign of $f_i$ on $\Delta$, and
$I^\Delta=\{i, \; \overline \Delta \cap H_i\neq \emptyset\}$ the
subset of hyperplans tangent to the domain $\Delta$. We set for
$\eta>0$
$$
\Delta^\eta=\{z\in \Delta, \;\; f_i(z)\epsilon_i^\Delta\ge \eta
\; \forall i\in I^\Delta\}.
$$
Since the integrant is exponentially decreasing at infinity, we
just have to evaluate the following integral
\begin{eqnarray*}
\vert\int_{\partial \Delta^\eta}\left( e^{-\lambda f_0}\prod \vert
f_i\vert^{\alpha_i}\right)\w_U\vert \le \sum_{i\in
I^\Delta}\vert\int_{\partial \Delta^\eta\cap
\{f_i\epsilon_i^\Delta=\eta\}}\left( e^{-\lambda f_0}\prod \vert
f_i\vert^{\alpha_i}\right)\w_U\vert .
\end{eqnarray*}
Now, if $i\in U$ then $\w_U$ vanishes on the set $\{f_i=\eta\}$.
On the other hand, if $i\not\in U$, then
$$
\int_{\partial \Delta^\eta\cap
\{f_i\epsilon_i^\Delta=\eta\}}\left( e^{-\lambda f_0}\prod \vert
f_i\vert^{\alpha_i}\right)\w_U\ \sim \eta^{\alpha_i}
\int_{\partial \Delta\cap H_i}\left( e^{-\lambda f_0}\prod \vert
f_i\vert^{\alpha_i}\right)\w_U,
$$
when $\eta$ tends to 0. Since the integral on $\partial \Delta
\cap H_i$ is finite since the weights $\alpha_i$ are all strictly
positive, we see that taking the limit $\eta\to 0$ we get by
Stockes theorem
$$
\lambda \int_\Delta \left( e^{-\lambda f_0}\prod_{j} \vert
f_j\vert^{\alpha_j}\right) df_0\wedge \w_U = \sum_{j\in U^c}
\epsilon(j,U)\alpha_j I_{\Delta, U\cup\{j\}}.
$$
We are now in a position to prove the result.
\begin{eqnarray}\nonumber
{dI_{\Delta,X}\slash d\lambda}&=&-\int_\Delta  \left( e^{-\lambda
f_0}\prod_{j} \vert f_j\vert^{\alpha_j}\right) f_0 \w_X
\\
\label{sum2}
 &=& -f_0(X) I_{\Delta, X} -\left( \sum_{j\in X}
c_{0,j} \int_\Delta  \left( e^{-\lambda f_0}\prod_{r} \vert
f_r\vert^{\alpha_r}\right) f_j \w_X\right).
\end{eqnarray}
Since $df_0=\sum_{j\in X} c_{0,j} df_j$, we see that
\begin{eqnarray*}
df_0\wedge \w_{X\setminus \ {j}}=c_{0,j} df_j\wedge
\w_{X\setminus\{ j\}}= \epsilon(j, X\setminus \{ j\}) c_{0,j} f_j
\w_X.\end{eqnarray*}
Hence, the sum in (\ref{sum2}) becomes
$$
 \sum_{j\in X} \epsilon (j,X\setminus\{j\}) \int_\Delta \left(
e^{-\lambda f_0}\prod_{r} \vert f_r\vert^{\alpha_r}\right) df_0
\wedge \w_{X\setminus \{j\}}.
$$
Using the cohomological relation we get
\begin{eqnarray*}
&& {dI_{\Delta,X}\slash d\lambda}
\\
&=& -f_0(X) I_{\Delta, X} -{1\slash \lambda} \left( \sum_{j\in X}
\sum_{r\in X\setminus \{j\}} \epsilon(j, X\setminus
\{ r\})\epsilon(r,X\setminus \{j\}) \alpha_r I_{\Delta,
X\setminus\{j\}\cup\{r\}} \right)
\end{eqnarray*}

 \ali\ali
Proof of proposition \ref{asymptotic}: A point $z$ in
$\ccc^\rho_X$ has the form
$$
z=X+\rho \sum_{j\in X} a_j e_{X\setminus \{j\}}, \;\;\;
(a_j)_{j\in X}\in (\R_+^*)^X.
$$
Thus, for $j\in X$ we have
$$
f_j(z)=\rho a_j l_j(e_{X\setminus \{j\}}),
$$
where $l_j$ is the linear form associated with $f_j$. For
$r\not\in X$
$$
f_r(z)= f_r(X)+\rho\sum_{j\in X} a_j l_r(e_{X\setminus \{j\}}),
$$
and
$$
f_0(z)= f_0(X) +\rho\sum_{j\in X} a_j,
$$
since by convention $l_0(e_{X\setminus \{j\}})=1$. Changing to the
variable $u_i=\lambda \rho a_j$ we see that if we set
$$
J_{X,X'}=\vert \det\left( l_r(e_{X\setminus\{j\}})_{{j\in X,\atop
r\in X'}}\right)\vert
$$
we get for all $\rho$ such that $\im(\rho)<0$ and $\lambda\not\in
i\R_-$
\begin{eqnarray*}
&&I_{X,X'}^+(\lambda)
\\
 &=&
 J_{X,X'} \lambda^{-\sum_{j\in X}(\alpha_j-\indic_{j\in X'}+1)}
 \left( \prod_{j\in X} \vert l_j(e_{X\setminus\{j\}})\vert^{\alpha_j-\indic_{j\in X'}} \right)
\left( \prod_{r\not\in X} \vert f_r(X)\vert^{\alpha_r-\indic_{r\in
X'}} \right)
\\
&& \;\;\;
  \int_{(\lambda \rho \R_+^*)^X}\left(
\prod_{j\in X} e^{-u_j} u_j^{\alpha_j-\indic_{j\in
X'}}\right)\left( \prod_{r\not\in X} h_r^{\alpha_r-\indic_{r\in
X'}}\right) \prod_{j\in X} du_j
\end{eqnarray*}
where
$$
h_r=1+\lambda^{-1} \sum_{j\in X} u_j
{l_r(e_{X\setminus\{j\}})\slash f_r(X)}
$$
(In $h_r^{\alpha_r}$ the determination of the logarithm is just
obtained by analytic extension of the logarithm, since at $u=0$,
$h_r=1$). Now, when $\lambda$ tends to infinity, then $h_r$
converges pointwise to 1. Using the dominated convergence theorem
we see that $I_{X,X'}(\lambda)$ is equivalent to
$$
D_{X,X'}e^{-\lambda f_0(X)}\lambda^{-\sum_{j\in X}
(\alpha_j-\indic_{j\in X'})}
$$
(and it can be made uniform in $\lambda$ in the domains $S^+$)
where $D_{X,X'}$ is the following constant
$$
J_{X,X'} \left( \prod_{j\in X}\Gamma(\alpha_j+\indic_{j\not\in
X'}) \vert l_j(e_{X\setminus\{j\}})\vert^{\alpha_j-\indic_{j\in
X'}} \right) \left( \prod_{r\not\in X} \vert
f_r(X)\vert^{\alpha_r-\indic_{r\in X'}} \right).
$$
Clearly, the term obtained for $X=X'$ is dominating and we get
that
$$
I_X(\lambda)\sim D_{X,X}\lambda^{-\alpha_X} e^{-\lambda f_0(X)}.
$$
where $D_{X,X}$ is as in proposition \ref{asymptotic}, since
$J_{X,X}=\prod_{j\in X} \vert l_j(e_{X\setminus \{j\}})\vert$.

\section{Acknowledgements}

The authors wish to thank A.N.Varchenko for helpful discussions.

\end{document}